\documentclass[a4paper, 11pt]{amsart}
\usepackage[latin1]{inputenc}
\usepackage[dvips]{graphicx}
\usepackage{amssymb}
\usepackage{amsmath}
\usepackage{amsthm}
\usepackage{amscd}
\usepackage{placeins}
\usepackage{xypic}

\theoremstyle{plain}

\newtheorem{proposition}{Proposition}[section]

\theoremstyle{definition}

\theoremstyle{remark}
\newtheorem*{remark}{Remark}
\newtheorem*{example}{Example}

\begin{document}
\title{Khovanov's conjecture over $\mathbb{Z}[c]$}
\author{Magnus Jacobsson}

\begin{abstract}
We disprove the conjecture by M.Khovanov [K1] on the functoriality 
of his link homology with polynomial coefficients. This is in contrast 
to the case of integer coefficients, where the functoriality was proven in [J].
\end{abstract}
\maketitle

\section{Introduction}
\subsection{Khovanov's Homology and Conjecture}
In his 1999 paper [K1] in Duke Mathematical Journal, Mikhail Khovanov
showed that the Jones link polynomial is 
the graded Euler characteristic of a bigraded homology module
$H^{i,j}(D)$ over $\mathbb{Z}[c]$, associated to a diagram $D$ of the link. 
$H^{i,j}(D)$ is the homology of a bigraded chain complex $C^{i,j}(D)$
with a $(1,0)$-bigraded differential.
He also explained how each link cobordism induces a homomorphism 
between the homology modules of its boundary links, and conjectured
that this homomorphism would be invariant up to sign 
under ambient isotopy of the link cobordism. 

In [J] we proved that, if formulated precisely, this conjecture is
indeed true, under the condition that the indeterminate $c = 0$, i.e. 
when the modules are abelian groups. (This case suffices to retrieve 
the Jones polynomial as the Euler characteristic.) In [K2], Khovanov
gave a different proof of this fact, and, in the introduction, renewed his 
conjecture in the general case.

\subsection{The Results of This Paper}
In the present note we consider the complete conjecture, with no
apriori conditions on $c$. We disprove it by giving an example
(Section \ref{proof}, Figure \ref{count}) of a trivial link cobordism 
from the unlink of two components to itself, which does not 
induce $\pm id$ on homology.

\subsection{A Remark on Coefficients}
Throughout this paper, we use coefficients in the polynomial ring (in
$c$) over $\mathbb{Z}_2 = \mathbb{Z} / {2\mathbb{Z}}$. This only
strengthens the result, and simplifies some calculations.

\subsection{Outline}
Section \ref{prel} contains necessary preliminaries on Khovanov homology.
In Section \ref{proof} we present the above-mentioned counterexample
to Khovanov's conjecture. 

\tableofcontents

\section{Preliminaries}
\label{prel}
\subsection{Khovanov's Chain Complex}

We briefly review the definition of Khovanov's chain complex
$C^{i,j}(D)$ and its differential. More details may be found 
in [V] or [J] (and, for the original definitions, in [K1]).

\subsubsection{The Chain Complex} 
\begin{figure}[htb]
\begin{center}
\includegraphics[width = 5 cm, height = 2 cm]{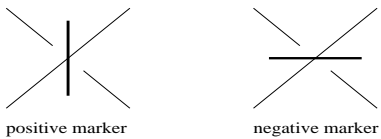}
\caption{Positive and negative markers.}
\label{markers}
\end{center}
\end{figure}

Let $D$ be a diagram of an oriented link. Recall that a (Kauffman) state of $D$ is 
a distribution of positive or negative {\em markers} to the crossings 
of the diagram (see Figure \ref{markers}). 
Each marker specifies one of the two ways to smooth the crossing. 
Smoothing all the crossings according to the markers at them gives a
set of embedded circles in the plane called the resolution of $D$
according to $S$.

Following Viro [V] we define an {\em enhanced state} $S$ to be a state
of $D$, together with an assignment of the symbols $\boldsymbol{1}$ or 
$X$ to each circle of the resolution of $D$ according to $S$. 
(In [V], $\boldsymbol{1}$ was denoted by the sign $-$ and $X$ by the
sign $+$.)

\begin{remark}
Since we will only be concerned with enhanced states, we call them states
from now on.
\end{remark}

Denote by $C(D)$ the free $\mathbb{Z}_2[c]$-module generated by all 
(enhanced) states of $D$. Denote by $w(D)$ the writhe of the diagram. 
Let $\sigma(S)$ be the sum of all signs of markers in the state $S$ 
and let $\tau(S)$ be the number of $X$:s minus the number of
$\boldsymbol{1}$:s on the resolution of $S$. 

$C(D)$ becomes a bigraded module $C = C^{i,j}(D)$, if we define the
grading parameters for an element $c^k S$ as 

\begin{displaymath}
\begin{split}
i(c^k S) &= \frac{w(D) - \sigma(S)}{2}\\
j(c^k S) &= - \frac{\sigma(S) + 2\tau(S) -3w(D)}{2} + 2k\\
\end{split}
\end{displaymath}
Notice that multiplication by $c$ affects only the second grading
parameter and that $\deg(c)= 2$.

\subsubsection{The Differential}
The differential $dS$ on a state $S$ in $C^{i,j}(D)$ is the sum of all
states $T$ in $C^{i+1,j}(D)$, which are {\em incident} to $S$ in the
following sense. 

\begin{figure}[htb]
\begin{center}
\includegraphics[width = 7 cm, height =  7 cm]{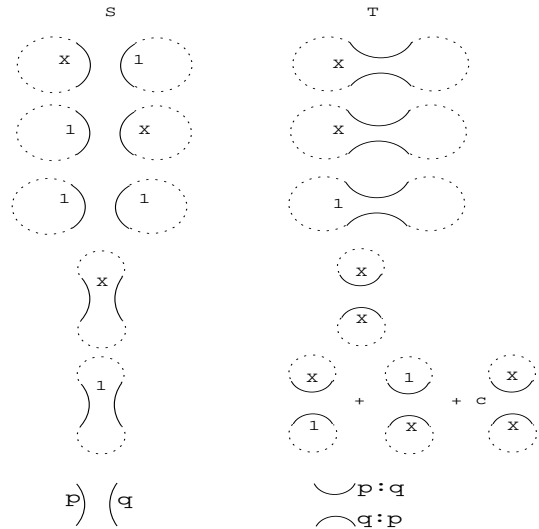}
\caption{Incident states}
\label{incident}
\end{center}
\end{figure}

\begin{itemize}
\item The markers of $T$ coincide with those of $S$ at all crossings
  except one, where $T$ has a negative marker and $S$ has a positive
  marker. This means that the resolution of $T$ has either one circle more or
  one circle less than the resolution of $S$.

\item The circles of the resolution which are common to $S$ and $T$
  are coulored with the same symbols in $S$ as in $T$. The circles which are not,
  have symbols related by the table in Figure \ref{incident}. It
  displays resolutions of $S$ and $T$ in a neighbourhood
  of the crossing at which the markers differ, with an indication
  (dotted arcs) of how the arcs are connected outside this
  neighbourhood. Note that the fifth row contains three different
  states incident to $S$. They appear in the differential as a
  sum. The last row is a mere mnemonic, summing up the other rows.

\item If $T$ is incident to $S$ then $c^k T$ is incident to $c^k S$
  for all $k$. 
\end{itemize}

\FloatBarrier
\subsection{The Second Reidemeister Move $\Omega_2$}
\label{basis}
In his proof of invariance of the (isomorphism class of) homology groups under 
Reidemeister moves, Khovanov shows the following proposition. Let $D'$ be obtained
from $D$ by a Reidemeister move of type 2 increasing the number of crossings.
(Reidemeister moves of type 2 will henceforth be referred to as
$\Omega_2$-moves).

\begin{proposition} There is a splitting $C(D') = C \oplus
 C_{contr}$, where $C_{contr}$ is chain contractible, and an
 isomorphism $\psi: C \xrightarrow{\cong} C(D)$. 
\end{proposition}

\begin{figure}[htb]
\begin{center}
\includegraphics[width = 7 cm, height = 2.5 cm]{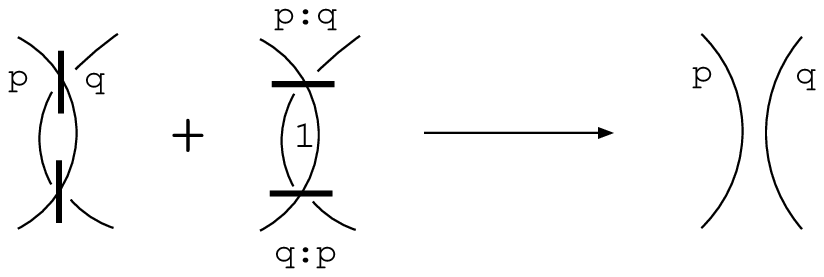}
\caption{The isomorphism $\psi$. The
  expressions on the left generate $C$.}
\label{psi2}
\end{center}
\end{figure}

\begin{figure}[htb]
\begin{center}
\includegraphics[width = 6 cm, height = 2.5 cm]{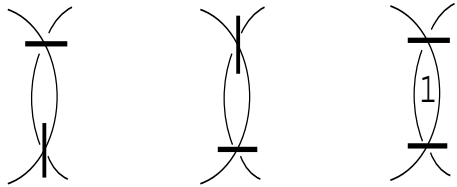}
\caption{Generators of the contractible splitting factor $C_{contr}$
 for the $\Omega_2$-move.}
\label{splitting2}
\end{center}
\end{figure}

Let us describe this splitting. Generators of $C$ and the isomorphism 
$C \cong C(D)$ are displayed in Figure \ref{psi2}. $C_{contr}$ is
generated (as a chain complex) by states as in Figure \ref{splitting2}.

\begin{remark} The symbols $p,q, p:q, q:p$ in the figures mean that
  the symbols on the arcs are related as in Figure \ref{incident}.
\end{remark}

\FloatBarrier

It follows that the composition $\Psi = \psi \circ pr_{C}$ given by 
\begin{displaymath}
C(D') = C \oplus C_{contr} \xrightarrow{pr_{C}} C \xrightarrow{\psi} C(D)
\end{displaymath}
is a chain equivalence. $pr_C$ denotes the canonical projection.
The inclusion $i_C$ and the inverse of $\psi$ give the homotopy inverse 
$\Psi^{inv}$. Using this it is straightforward to prove (using the
same method as in Subsection 2.4 in [J]) the following proposition.

\begin{figure}[htb]
\begin{center}
\includegraphics[width = 6.5 cm, height = 2 cm]{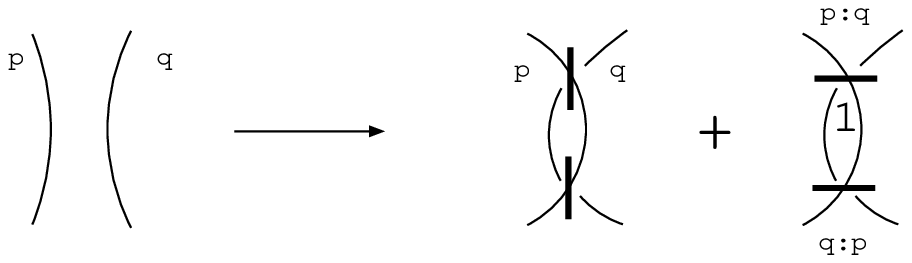}
\caption{The chain equivalence $\Psi^{inv}$.}
\label{chain2c}
\end{center}
\end{figure}

\begin{figure}[htb]
\begin{center}
\includegraphics[width = 6.5 cm, height = 8 cm]{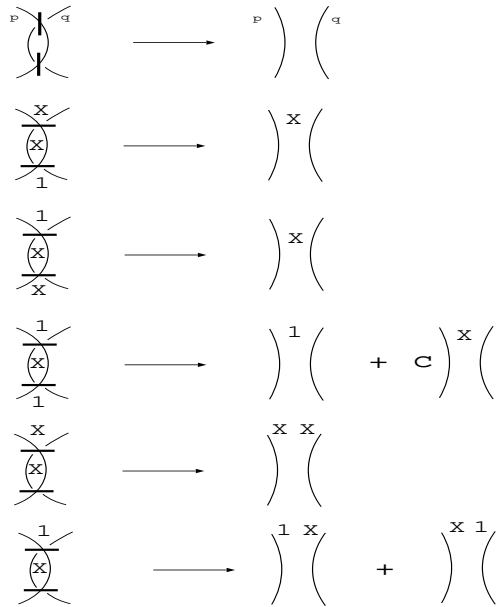}
\caption{The effect of the second Reidemeister move on states. 
The chain equivalence $\Psi$. A state with any other local
configuration is mapped to zero.}
\label{chain2c22}
\end{center}
\end{figure}

\begin{proposition}
\label{prop}
The chain equivalences $\Psi$ and $\Psi^{inv}$  above have the values
given in Figures
\ref{chain2c22} and \ref{chain2c} on states with given local configuration. 
\end{proposition}

\FloatBarrier
\subsection{Khovanov's Conjecture}

Recall that a {\em link cobordism} $\Sigma \subset \mathbb{R}^3 \times
[0,1]$ between links $L \subset \mathbb{R}^3 \times \{0\}$ 
and $L' \subset \mathbb{R}^3 \times \{1\}$ is a smoothly embedded compact oriented 
surface whose boundary consists of the links.

Every link cobordism can be presented by a {\em movie}, which is the
sequence of diagrams of the intersection of $\Sigma$ with the constant 
time hyperplanes $\mathbb{R}^3 \times \{t\}$ (subject to certain
genericity assumptions). For all but a finite number of values of $t$ 
(critical levels), this intersection is a link diagram, and at times 
just before and after a critical level the diagrams differ by a 
Reidemeister move or a Morse modification. Between two critical levels
the diagram experiences a planar isotopy.

When a link cobordism is altered by an ambient isotopy, its movie changes accordingly.
(For additional details on link cobordisms, see e.g. [J].)

In his paper [K1], Khovanov described chain equivalences induced by all 
Reidemeister moves and chain maps induced by Morse modifications on
the link diagram. (In the previous section we explained how this is done for the
second Reidemeister move, since this is all we will need in this paper.)
It is also clear from the construction that states can be traced
through any planar isotopy between two diagrams an hence induces a
chain isomorphism between the associated chain complexes.

It follows that each movie presentation of a link
cobordism induces a map between the homologies of its boundary link
diagrams, via composition.

Khovanov conjectured that this map be invariant up to sign under 
ambient isotopy of the link cobordism. We pointed out in [J] that 
it is necessary to require stability of the boundary during 
the isotopy. (For example, that the boundary is held fixed during the
isotopy.) We then proved the conjecture under the assumption 
that $c=0$.

In the next section we prove that the assumption $c = 0$ is in fact 
necessary, by exhibiting a movie presentation of a trivial link
cobordism which does not induce the identity on homology.
\FloatBarrier

\section{A Counterexample to Khovanov's Conjecture}
\label{proof}

\begin{figure}[htb]
\begin{center}
\includegraphics[width = 10 cm, height = 4 cm]{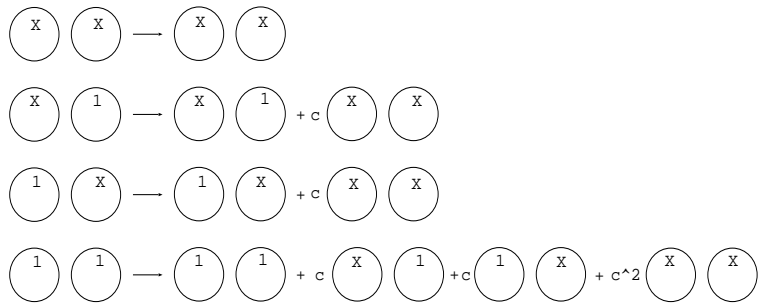}
\label{count2}
\caption{The induced map on homology, forcing $c$ to be $0$.}
\end{center}
\end{figure}

\begin{example}
\label{example}
To begin with, let us compute the effect on chains under the following
sequence (*) of $\Omega_2$-moves. The sequence starts with a crossing 
somewhere in a link diagram. Just above the crossing an $\Omega_2$-move
is performed, creating two bigons. Removal of the lower one of these
via an $\Omega_2$-move in the opposite direction returns us to the original diagram:

\begin{figure}[!h]
\begin{center}
\includegraphics[width = 5 cm, height = 2.5 cm]{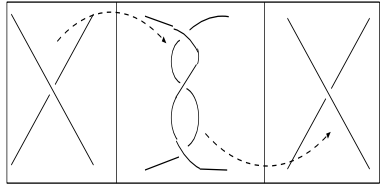}
\caption{The movie (*).}
\label{der}
\end{center}
\end{figure}

Figure $\ref{map14}$ shows what happens on $C^{i,j}(D)$. We include
only states whose images we will need later. To verify
this table is straightforward, given Proposition \ref{prop}. 

\begin{figure}[htb]
\begin{center}
\includegraphics[width = 13 cm, height = 8 cm]{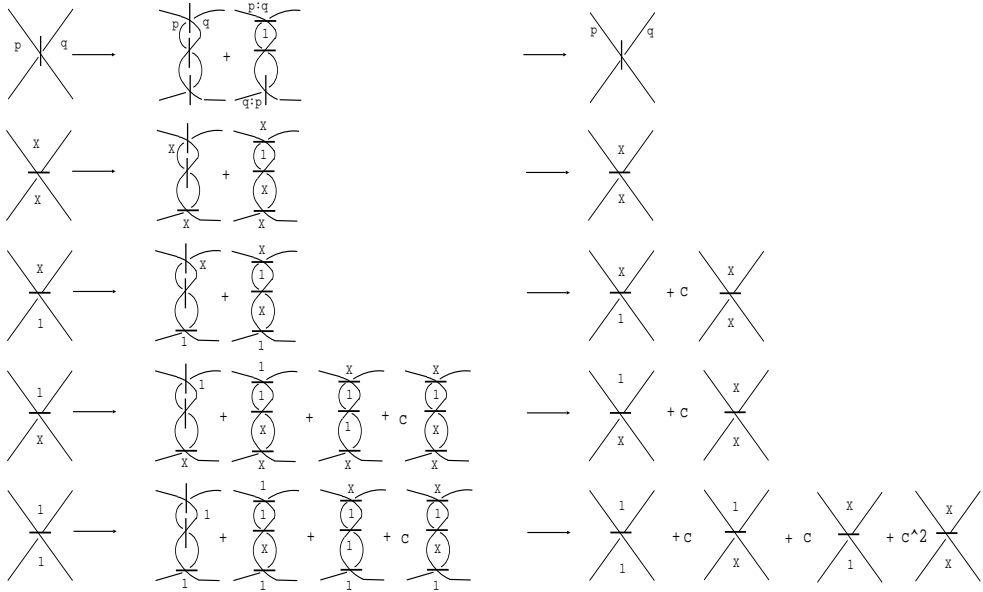}
\caption{The chain map induced by (*).}
\label{map14}
\end{center}
\end{figure}

Now consider the sequence of moves in Figure \ref{count}. 
It starts and ends with the trivial unlink diagram of two components. 
Since evidently the differential is zero for this diagram, the 
homology modules are canonically isomorphic to the chain modules.
Each state is just a distribution of $\boldsymbol{1}/X$:s to the components.

The sequence starts with an $\Omega_2$-move on the trivial diagram, 
sliding the left circle above the right one. This induces the chain equivalence 
$\Psi^{inv}$ as described in Proposition \ref{prop}. Then the sequence (*)
described above is applied at the upper crossing. Finally the circles
slide back again inducing $\Psi$. 

This is a movie of a link cobordism from the trivial unlink to itself, 
which is isotopic to the cylinder on the unlink. A simple computation 
shows that the map induced on the homology of the unlink is the one 
given in Figure \ref{count2}. Clearly it is only $id$ if $c = 0$.
This disproves the general version of Khovanov's conjecture.
\end{example}

\begin{figure}[htb]
\begin{center}
\includegraphics[width = 5 cm, height = 4 cm]{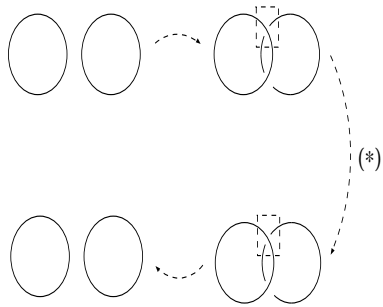}
\caption{A movie of a trivial link cobordism. (*) refers to the movie
  in Figure \ref{der} applied in the dotted square.}
\label{count}
\end{center}
\end{figure}

\FloatBarrier

\end{document}